\documentclass[11pt]{article}

\usepackage{amsfonts,amssymb}

\newcommand{\Z}{\mathbb{Z}}
\newcommand{\R}{\mathbb{R}}

\input{epsf.sty}

\begin{document}

\title{A geometric method to compute some elementary integrals}
\author{J. Scott Carter \\ Abhijit Champanerkar \\University of South Alabama }

\maketitle

\begin{abstract}
An elementary, albeit higher dimensional, argument is used to compute the area under the power function curve between $0$ and $1$.
\end{abstract}

\section{Introduction}

The purpose of this paper  is to give a geometric proof that $\int^1_0 x^n \ dx = \frac{1}{n+1}$ for $n\ge 0$. The proof, although quite simple, depends deeply on the geometry of higher dimensional space.
The argument is designed to inform the reader, presumably an advanced undergraduate, about some of the geometrical configurations that can be found in higher dimensional space. Thus, our point of view is to take a familiar concept --- that of area under a curve --- and reveal it in a different light. 

In order to analyze higher dimensional geometry, we need to use rigid transformations of these spaces. Such transformations are found among the set of square matrices. Thus our proof uses basic linear algebra to check all the details. It is our hope that by thinking geometrically about the matrix manipulations, the reader can gain geometric intuition about higher dimensional space.

Let us sketch the proof. Let $x \in [0,1]$, then the quantity $y=x^n$ is the $n$-dimensional volume of an $n$-cube each of whose edge lengths is $x$. Such an $n$-cube can be found embedded in the unit $(n+1)$-cube, at height $x$ with its edges parallel to the coordinate vectors $e_1=(1,0,\ldots,0)$, $e_2=(0,1,0,\ldots, 0)$, through $e_n=(0,0,\ldots, 1,0)$. That is to say that for $x$ between $0$ and $1$, the cube rests in the hyperplane $x_{n+1}=x$ in a standard position. As $x$ varies from $0$ to $1$, the union of these $n$-cubes forms a pyramid with $n$-cubical base. The area under the curve $x^n$ coincides with the $(n+1)$-dimensional volume of the pyramidal cone on the unit base. Our main task then will be to show that $(n+1)$-of these pyramids fill the unit $(n+1)$-cube.
To see this, we will show that the pyramid intersects a simplex dual to the diagonal line from
$(0,0,\ldots,0,0)$ to $(1,1,\ldots , 1, 1)$ in the barycenters of some of its lower dimensional faces. The simplex has a cyclic symmetry of degree $(n+1)$-and the rotations can be applied to the pyramid
so that its orbits fill the $(n+1)$-cube.
Specifically, there is a symmetry of the unit $(n+1)$-cube of order $n+1$ which gives
an action of the cyclic group, $\Z_{n+1}$, on the unit $(n+1)$-cube. The orbit of the pyramid under this action consists of (n+1) copies of the pyramid, which tesellate the $(n+1)$-cube. Thus the pyramid takes up $1/(n+1)$ of the volume thereof. Since the volume of the pyramid coincides (up to choice of unit) with the integral, $\int^1_0 x^n \ dx$, this gives the result.

In the preceding paragraph, a lot of terminology was introduced without explanation. Before we embark on the full higher dimensional journey, we first run through the idea in dimension 3. Then we examine the ideas in $2$-dimensions. After those examples are fully worked, we introduce the general terminology in higher dimensional space, and describe the intersections between certain subspaces in higher dimensional space. Then we describe how the pyramidal shape can be rotated to fill the cube. We discuss filling hyper-cubes with simplices, we give some illustrations, and we give the {\sc Mathematica} code that can be used to generate pictures of the pyramidal sets.

\subsection*{Acknowledgments}

Seiichi Kamada gave us  the description in terms of the link of the vertex, and Makato Matsumoto gave us the description of the decomposition of the $(n+1)$-cube as a union of $(n+1)!$ simplices.  Susan Williams gave us the anyaltic description of Section~\ref{last}. JSC is supported by NSF grant 0301095; AC is supported by NSF grant DMS-0455978.

\section{Three dimensions}
\label{3d}

Here is the second most easy case. The quantity $y=x^2$ represents the area of a square all of whose edge lengths are $x$.
We arrange a set of these squares in a unit cube as follows. The unit square lies on the top of a $1 \times 1 \times 1 $ cube as the set $S_1=\{ (u,v,1) : 0\le u \le 1 \  \ \& \ \ 0 \le v \le 1 \}$.
The square of side length $x$, and area $x^2$ is found as $S_x=\{(u,v,x):
 0\le u \le x \  \ \& \ \ 0 \le v \le x \}.$ The union of these squares $\cup_{0 \le x \le 1} S_x$
forms a pyramid with a square base, height $1$ and volume $\int_0^1 x^2 dx$.
The pyramid has two faces that are right isosceles triangles with sides of length $1$, $1$ and $\sqrt{2}$. 
For either of the remaining two right triangular faces, the hypotenuse is 
of length $\sqrt{3}$, a side of length $1$, and one of length
$\sqrt{2}$. The pyramid that we are considering, then is the intersection of all the convex sets in $3$-space that contain the vertices $(0,0,0)$, $(1,0,1)$, $(0,1,1)$, $(1,1,1)$, and $(0,0,1)$.

\begin{figure}[htb]
\begin{center}
\mbox{
\epsfxsize=3in
\epsfbox{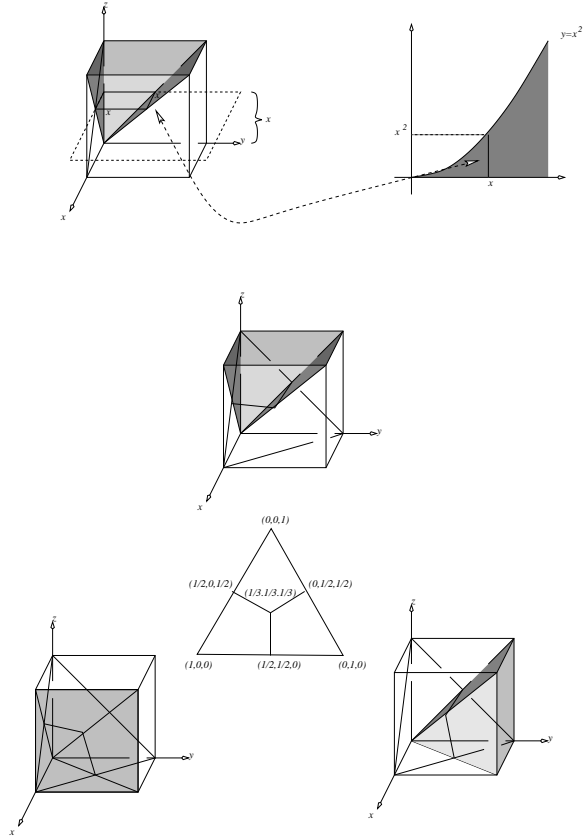}
}
\end{center}
\caption{Tiling the three dimensional cube with pyramids}
\label{threecubes}
\end{figure}

We call the edge from $(0,0,0)$ to $(1,1,1)$ the {\it main diagonal}. It can be given in parametric form as $\ell(t) = (t,t,t)$ for $t\in [0,1]$. The main diagonal intersects an equilateral triangle in the cube at the centroid of the triangle.

The equilateral triangle in question has as its vertices $(1,0,0)$, $(0,1,0)$, and $(0,0,1)$.
Observe that the intersection of the triangle with the diagonals of the pyramid of length $\sqrt{2}$ occur at the points $(1/2,0, 1/2)$, and $(0,1/2,1/2)$ while  the big diagonal intersects at $(1/3,1/3,1/3)$. The triangle can be given as the intersection of the cube with the plane $x+y+z=1$. Each of  the little diagonals can be parameterized, and these intersections can be determined. The entire pyramid intersect the triangle in a kite shape with vertices $(0,0,1)$, $(1/2,0,1/2)$, $(0,1/2,1/2)$, and $(1/3,1/3,1/3)$. For convenience, we indicate the kite as $D_2[3]$ and write $$D_2[3]=H(\{(0,0,1),(1/2,0,1/2), (0,1/2,1/2),(1/3,1/3,1/3)\});$$ the subscript $2$ indicates that the kite is two dimensional, and the $[3]$ indicates that the third component of each vertex is non-zero. The letter $H$ indicates the convex hull which will be defined in general in Section~\ref{nd}.

Under the rotation 
$$(1,0,0) \mapsto (0,1,0) \mapsto (0,0,1) \mapsto (1,0,0),$$ 
the kite $D_2[3]$ rotates to $D_2[1]=H(\{(1,0,0),(1/2,1/2,0),(1/2,0,1/2),(1/3,1/3,1/3)\})$, and this in turn rotates to
$D_2[2]=H(\{(0,1,0),(0,1/2,1/2),(1/2,1/2,0),(1/3,1/3,1/3)\}).$ These three kites tile the equilateral triangle. 
The restriction of the rotation to the pyramid gives two more pyramids. The three of these tile the cube.
Under the rotation, the square base rotates first to the $x=1$ face and then to the  $y=1$ face of the cube. The $y=0$ face of the cube is tiled by two isosceles right triangles as are the faces $x=0,$ and $z=1$. The illustration in Fig.~\ref{threecubes} gives all the details.

What we have shown then is that the integral $\int^1_0 x^2 \ dx = 1/3$. The quantity on the left side represents the volume of the pyramid, and the quantity on the right represents $1/3$ of the cube.

\section{The two dimensional case}
\label{2d}

\begin{figure}[htb]
\begin{center}
\mbox{
\epsfxsize=2.5in
\epsfbox{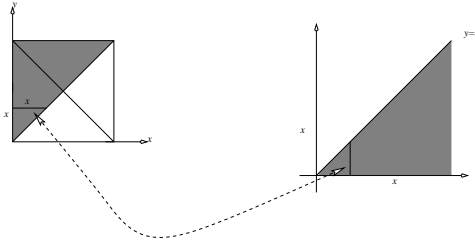}
}
\end{center}
\caption{Tiling the square with triangles}
\label{twocubes}
\end{figure}

The first sentence of the preceeding section was a teaser. If the case $n=2$ is the second most easy case, what is the easiest case? For $n=1$, we arrange horizontal line segments of length $x$ at height $y=x$ in the plane. The union of these segments is the isosceles right triangle with vertices $(0,0)$,  $(0,1)$, and $(1,1)$. The hypotenuse of this triangle is the line segment $\ell(t)=(t, t)$ for $t\in [0,1].$ The triangle clearly occupies half of the square, $\{(x,y) \in \R^2  :  0\le x \le 1 \ \  \& \  \ 0 \le y \le 1 \}$. Moreover, we can reach the other half of the triangle  by reflecting in the line $\ell(t)$. Fig~\ref{twocubes} contains the details of the construction.
The analogue of the triangle is the line segment on $x+y=1$ that lies in the first quadrant (so $0 \ge x,y$). 

We remark here on an interesting difference between dimensions $2$ and $3$. In dimension $2$ the symmetry of the square that takes the lower triangle to the upper one is a reflection in the diagonal. In dimension $3$ it is a rotation. If we extend the former situation to all of $3$-dimensions, then
we can think of the reflection in the line as having been achieved by a rotation in $3$-space. In the general situation,  the even dimensional cubes will be reflected.

\section{General notation}
\label{nd}

\subsection{The $(n+1)$-dimensional cube}

We begin by describing the unit $(n+1)$-dimensional cube. This is given abstractly as the Cartesian product: $$\underbrace{[0,1] \times[0,1] \times \cdots \times [0,1]}_{(n+1)-{\mbox{\rm factors}}}= [0,1]^{n+1}$$
or more concretely as a subset of $\R^{n+1}$
$$C^{n+1}= \{ (x_1,x_2,\ldots, x_{n+1}) \in \R^{n+1}: \ 0\le x_j \le 1 \ \ {\mbox{\rm for all }} \ \  j=1,2, \ldots, n+1 \}.$$
In $\R^{n+1}$ it is  convenient to denote the standard basis vectors as follows:
$e_1=(1,0,\ldots ,0),$ $e_2=(0,1,\ldots ,0 )$, $\ldots,$ $e_n=(0,0,\ldots, 1, 0)$, and $e_{n+1}=(0,0,\ldots, 0, 1)$. In general, throughout our  discussion $n$ is a fixed integer and $n\ge0$. In the case that $n=0$, the standard basis vector is $1$.

The important properties of a {\it  basis} such as  ${\mathcal E} = \{ e_1, e_2, \ldots, e_{n+1} \}$ are
\begin{enumerate}
\item
Any vector $v \in \R^{n+1}$ can be written as a linear combination --- $v= \sum_{j=1}^{n+1} \alpha_j e_j$ where $\alpha_j \in \R$ of elements in ${\mathcal E}$;
\item Such an expression is unique. That is, 
if $v = \sum_{j=1}^{n+1} \alpha_j e_j = \sum_{j=1}^{n+1} \beta_j e_j,$ then
$\alpha_j =\beta_j$ for all $j=1,2, \ldots, n+1$. 
\end{enumerate}

A set, $S$,  is {\it convex}  if any two points that are in the set can be joined by a line segment that is contained in $S$. If $S\subset \R^{n+1}$, and $P,Q \in S$, then the line segment from $P$ to  $Q$ can be given in parametric form  as $\ell_{PQ}(t)= P + t(Q-P)$, for $t\in [0,1]$. Here we are using the vector space structure of $\R^{n+1}$ so that if $P=(x_1, x_2 ,\ldots, x_{n+1})$ and $Q=(y_1,y_2,\ldots , y_{n+1})$, then $Q-P= (y_1-x_1, y_2-x_2, \ldots, y_{n+1}-x_{n+1})$. A point on the line segment $\{ \ell_{PQ}(t): t \in [0,1]\}$ is determined by vector addition as well. In the $j\/$th  coordinate, we have $x_j+t(y_j-x_j)=(1-t)x_j + t y_j$. At $t=0,$ the point $\ell_{PQ}(0)=P$ and  $\ell_{PQ}(1)=Q$. If $S$ is convex, then  $\ell_{PQ}(t) \in S$ for all $t\in [0,1]$.

The {\it convex hull} of  a set,  $V$, is the intersection of all the convex sets that contain $V$.
For example let $J'=\{j_1, j_2, \ldots, j_k\}$ denote any $k$-element subset of $\{1, 2, \ldots, n+1\}$. For  ease of notation, we assume  that $j_1 < j_2 < \ldots < j_k$. Then let
$e_{J'}= \sum_{i=1}^k e_{j_i}$. The vector $e_{J'}$ is the vector that has $1$ in each of  the entries $j_i$ and $0$ elsewhere. Of course, $e_{\{j\}}=e_j,$ and in this case we use the latter notation. By  convention, $e_{\emptyset}=(0,0,\ldots, 0)$. We consider,
$V = \{ e_{J'}: J' \subset \{1, 2, \ldots, n+1\} \}.$
The convex hull of $V$ is the unit cube.  Denote the convex hull of any set $S$, by $H(S)$.

The construction of the preceeding paragraph indicates the deep connections among the $2^{n+1}$ vertices of the $(n+1)$-cube, the $2^{n+1}$ subsets of $\{1,2, \ldots, n+1 \}$, and the lattice of inclusions among these subsets. Specifically, if $J^{\prime \prime}= J^{\prime}\cup \{i_{k+1}\}$, then there is an edge of length $1$ in the $(n+1)$-cube that  connects the corresponding points
$e_{J^{\prime \prime}}$ and $e_{J^{\prime}}$.

\subsection{The $(n+1)$-dimensional pyramid}
As a further example, consider any subset $J=\{j_1, j_2, \ldots, j_k \}\subset \{1,2,\ldots, n\}$, and let
$V=\{e_{J\cup \{n+1\}}:  J \subset \{1, 2, \ldots, n\} \}.$ Then $V$ consists of the union of the vertices of the $n$-cube, $C^{n}(1) = H(V).$  The $1$ is meant to denote that this cube is the ``ceiling'' of the $(n+1)$-cube. Now let ${\mbox{\rm Py}}_{n+1}$ denote the convex hull of $V\cup \{e_\emptyset \}$. The set ${\mbox{\rm Py}}_{n+1}$ is a ``pyramidal'' set that consists of the set of all lines segments from a point in the ``horizontal'' $n$-cube $C^n(1)=\{(x_1,x_2,\ldots, x_n,1)\in R^{n+1}: 0 \le x_j \le 1 \ \  {\mbox{\rm for all}} \ \  j= 1, 2, \ldots , n \}$ to the origin $(0,0,\ldots,0)$. The intersection of the pyramid ${\mbox{\rm Py}}_{n+1}$ with the horizontal
hyperplane $x_{n+1}= x$ is an $n$-cube with the length of each of its edges equal to $x$. Thus the $(n+1)$-dimensional volume of ${\mbox{\rm Py}}_{n+1}$ coincides with
$\int^1_0 x^n \ dx$.

Let us examine in some detail the claim that the intersection of the pyramid
${\mbox{\rm Py}}_{n+1}$ with the hyperplane $x_{n+1}=x$ is an $n$-cube with each of the lengths of its edges equal to $x$. Let $P=(x_1, x_2, \ldots, x_n, 1)$ denote a point in the horizontal $n$-cube, so that for $j=1, 2, \ldots , n$, we have $0\le x_j \le 1.$ The pyramid consists of all lines from the base cube to the vertex $Q=(0,0,\ldots, 0).$ Such a line would be
$$\ell_{PQ}(t)= (tx_1, tx_2, \ldots, tx_n, t).$$ If $t=x$, since $0<x_j<1$, we have $0< x x_j < x$, for $j=1, 2, \ldots , n$.  That is, each of the early coordinates is bounded above by $x$.

Now we examine the claim that the $(n+1)$-dimensional volume of  ${\mbox{\rm Py}}_{n+1}$ is
$\int^1_0 x^n \ dx$. For each $x \in [0,1]$ the height of the curve $x^n$ is the $n$-dimensional volume of the $n$-cube each of whose edge lengths is $x$. Cavalieri's principle (see any standard calculus text) states that if two figures have the same cross sectional measure then the integrals coincide. This idea was used by Archimedes, for  example, to show that the volume of  a sphere could be obtained from  the volume of a cylinder and the volume of a cone. Here we use this idea to translate the integral $\int^1_0 x^n \ dx$ into a computation of the volume of  the pyramid ${\mbox{\rm Py}}_{n+1}$. Upon careful consideration, it seems much more natural to consider this rectilinear figure rather than the curved region bounded above by the power function.  The power function is an expression of $n$-dimensional volume. The integral adds these volumes to  get  an $(n+1)$-dimensional volume.

\subsection{An $n$-dimensional simplex}

Among the pantheon of geometric figures in the plane, the square and the equilateral triangle stand out. The $3$-dimensional  analogue of a square is a cube, and the $3$-dimensional analogue of a triangle is a tetrahedron. The $n$-dimensional simplex, or $n$-simplex, is the higher dimensional analogue of the triangle and tetrahedron. It is a regular figure and can be thought of as the controlling figure of higher dimensional geometry in the following sense.

In Euclidean geometry, we learn that any two points determine a line, and any three non-collinear points determine a plane. In such a plane, there is a triangle  that is the convex hull of the three given points. Any four points that are not coplanar determine, via their convex hull, a tetrahedron. And that's as far as things go if you limit yourself to  exist in the $3$-dimensional world.  If on the other hand, you see that the analogy should carry on, we could consider any five points  that do not co-exist in a single $3$-dimensional hyperplane  should form the basic figure, a $4$-simplex, in $4$-space. And so forth.

In the triangle each pair of points is joined by an edge, and the same is true of the tetrahedron. It has six edges, four triangular faces, and four vertices. In an $n$-simplex every pair  of vertices is joined. Here is a nice way to draw the $1$-dimensional spine (edges) of the $n$-simplex.  Consider the set of $(n+1)$-points arranged evenly spaces on a unit circle. Then draw a line between every pair of points.  For example, the $4$-simplex projects to a regular pentagram (see Fig~\ref{petnagram}. It is amusing and rewarding to draw a circle with a protractor, measure an angle of $360/(n+1)$, mark points on the circle subtended by multiples of these angles and, draw lines between every  pair of points.  The observant artist will  see other lower dimensional simplices emerge from such a process.

\begin{figure}[htb]
\begin{center}
\mbox{
\epsfxsize=2.5in
\epsfbox{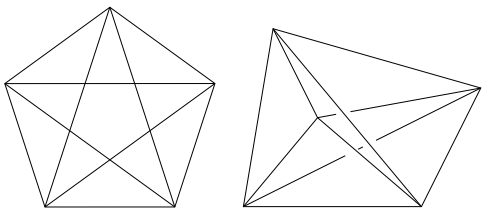}
}\end{center}
\caption{The four dimensional simplex from two standard points of view}
\label{pentagram}
\end{figure}

\begin{figure}[htb]
\begin{center}
\mbox{
\epsfxsize=4in
\epsfbox{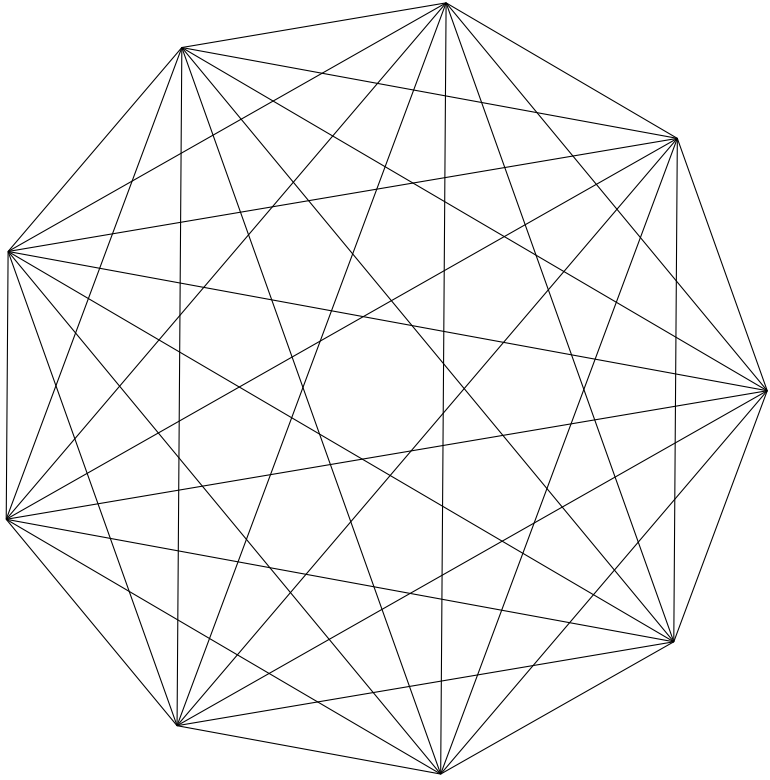}
}
\end{center}
\caption{The set of edges in the $8$-dimensional simplex}
\label{8simplex}
\end{figure}

An illustration of an $8$-simplex is given in Fig.~\ref{8simplex}.

We are interested in a particular regular $n$-simplex, that is the convex hull of a set  of points in the $(n+1)$-cube. The points in question consist of $e_1=(1,0, \ldots ,0, 0)$, $e_2=(0,1,\ldots, 0,0)$, through
$e_{n+1}=(0,0, \ldots, 0, 1)$.

In an $n$-simplex every pair of vertices are joined by a line segment. The distance from $e_i$ to $e_j$, computed using $d(e_i,e_j)= \sqrt{ (\Delta x_1)^2 + (\Delta x_2)^2 + \cdots + (\Delta x_{n+1})^2)} =\sqrt{2}.$
 Similarly, every set of $3$-vertices determines a triangle, every set of $4$-vertices determines a tetrahedron, and in general every set of $k$-vertices determines a $(k-1)$-dimensional simplex that is the convex hull of that subset. Each of these sub-figures is also regular. For definiteness, we label the $n$-simplex in which we are interested as the convex hull of $\{ e_1, e_2, \ldots, e_{n+1} \}$
$$S_n= H(\{ e_1, e_2, \ldots, e_{n+1} \}).$$

\section{The intersection between ${\mbox{\rm Py}}_{n+1}$ and $S_n$}

Before we make the computations of this section, we want to remind you of where we are in this process. We have argued that  $$\int^1_0 x^n \ dx ={\mbox{\rm Vol}}({\mbox{\rm Py}}_{n+1})$$ where by ${\mbox{\rm Vol}}$, we mean the $(n+1)$-dimensional volume of that figure. We want to compute this volume by showing that the pyramid takes up exactly, $1/(n+1)$ of the $(n+1)$-cube. To see this, we are going to compute the intersection of the pyramid with the simplex, and  apply a set of symmetries of the  simplex to the entire $(n+1)$-cube. As the symmetry is applied the pyramid will intersect another portion of the simplex, and this intersection is an indication of how the pyramid and its translations fill the cube. Now we commence.

Let $j$ denote a fixed integer from $1$ to  $n$. Consider the edge $\ell_j$ in the pyramid ${\mbox{\rm Py}}_{n+1}$ that consists of the line segment from $e_{\emptyset}$ to  $e_{j}+e_{n+1}$. In parametric form,
$\ell_j(t) = t (e_j + e_{n+1})$ for  $t\in [0,1]$. Meanwhile, the edge in the simplex that joins $e_j$ to $e_{n+1}$ is given as $S_j(u)= e_j + u(e_{n+1}-e_j)$ for  $u \in [0,1]$. The intersection between these two lines occurs when
$e_j + u(e_{n+1}-e_j)=t (e_j + e_{n+1})$. Since $\{ e_j: j=1, 2, \ldots ,  n+1 \}$ is a basis for the vector space $\R^{n+1}$,  we have $1-u=t$, and $u=t$. Therefore  the intersection point is $(e_j+e_{n+1})/2.$

Let $J=\{j_1 , j_2, \ldots, j_k \} \subset \{1,2,\ldots, n\}$. Then $e_{J\cup\{n+1 \} }$ is a vertex  of the $n$-dimensional cube $C^n(1)$ that is the ceiling of the $(n+1)$-cube. The line
$\ell_J(t) = t e_{J\cup\{n+1 \} } = t e_J + t e_{n+1}$ is the line on the cone from the origin to that vertex 
--- it is an edge of  the pyramid ${\mbox{\rm Py}}_{n+1}.$   This edge intersects the simplex $S_n$ in a $k$-dimensional  face  that is determined by the subset $J=\{ j_1, j_2, \ldots, j_k \}$. 
Let $$T_{J\cup \{n+1\}}= H(\{ e_{j_1} , e_{j_2},\ldots, e_{j_k}, e_{n+1} \}).$$
Then $T_{J\cup \{n+1\}}$ has the geometric shape of a $k$-simplex. For example, if $k=1$, then $T_{\{j_1,n+1\}}$ is  a line segment of length $\sqrt{2}$. If $k=2$, then $T_{\{j_1,j_2,n+1\}}$ is an equilateral triangle with vertices $e_{j_1}$, $e_{j_2}$
and $e_{n+1}$.
A linear combination, $\sum_{i=1}^k u_i e_{j_i} + u_{n+1}e_{n+1}$ is an element of $T_J$ provided $\sum_{i=1}^k u_i +u_{n+1}=1$ (and 
$u_i \geq 0$ ).

The intersection between $T_{J\cup\{n+1\}}$ and $\ell_J(t)$ is computed:
$$t\sum_{i=1}^{k} e_{j_i}+t e_{n+1} =\sum_{i=1}^k u_i e_{j_i} + u_{n+1}e_{n+1}.$$
We equate the coefficients of $e_{j_i}$ to get $u_i=t$ for each $i=1,2, \ldots, n+1$. Since 
$\sum_{i=1}^k u_i +u_{n+1}=1$, $t=1/(k+1)$ and the intersection point is
$$p_{J\cup\{n+1\}}=1/(k+1) e_J + 1/(k+1) e_{n+1}= 1/(k+1) \sum_{i=1}^k e_{j_i} + 1/(k+1) e_{n+1}.$$
The reader will observe that this computation used the fact that ${\mathcal  E}=\{e_1, e_2, \ldots, e_{n+1} \}$ is  a basis, so that there is a solution to these equations and the solution is unique.

For $J =\{ j_1, j_2 , \ldots, j_k\}\subset \{1,\ldots, n+1\}$, let $p_{J
} = \frac{1}{k+1} \sum_{j \in J} e_j$. The intersection between ${\mbox{\rm Py}}_{n+1},$ and $S_n$ consists of the set
$$D_n[n+1] = H(\{p_J: n+1 \in J \subset \{1,2,\ldots, n+1\} \}).$$
The set $D_n[n+1]$ is an $n$-dimensional cubiod with its two dimensional faces being kites. It is combinatorially equivalent to the $n$-cube, $C^n(1)$ on the ceiling.
More generally, let $D_n[i] = H(\{ p_J: i \in J \subset \{1,2, \ldots, n+1 \}\})$. Each $D_n[i]$ is an 
$n$-cuboid, and the union of these for $i=1, \ldots, n+1$ tiles $S_n.$

\section{Symmetry of the $(n+1)$-cube}

Recall that the line $\ell_{\{1, 2, \ldots, n\}}(t)=(t,t,\ldots , t)$  is the main diagonal of the $(n+1)$-cube. It  joins the 
points $(0,0,\ldots, 0, 0)$ and $(1,1, \ldots, 1, 1)$. We will rotate the cube about the axis 
$\ell_{\{1, 2, \ldots, n\}}(t)$ 
in such a way that $e_1$ rotates to $e_2$, $e_2$ rotates to $e_3$, and so forth, $e_{n+1}$ rotates back to $e_1$.
The rotation which acheives this is given by

$$\Theta_{n+1}= \left( \begin{array}{cccccc}
0 & 0 & \ldots & 0 & 0 & 1 \\
1 & 0 & \ldots & 0 & 0 & 0 \\
0 & 1 & \ldots & 0 & 0 & 0 \\
\vdots & \vdots &   & \vdots & \vdots & \vdots \\
0 & 0 & \ldots & 1 & 0 & 0 \\
0  & 0 & \ldots & 0 & 1& 0 \end{array}  \right) $$
With $\Theta_{n+1}$ given in this  form,  for  $n$ an even integer, $\det{(\Theta_{n+1})}=-1$, for $n$ an odd integer, $\det{(\Theta_{n+1})}=1$, and in general
$(\Theta_{n+1})^{n+1}={\mbox{\rm Id}}$. In other words 
$\Theta_{n+1}(x_1,x_2, \ldots ,x_n,x_{n+1})=(x_{n+1},x_1,x_2, \ldots ,x_n)$.

Note that $\Theta_{n+1}$ is a symmetry of the $(n+1)$-cube since it preserves the set of vertices. Since $S_n$ passes through the points $e_1$ to $e_{n+1}$, and $\Theta_{n+1}$ preserves the set $\{e_1,\dots,e_{n+1}\}$, 
$\Theta_{n+1}$ is a symmetry of $S_n$. 

We can see that the  action of $\Theta_{n+1}$ on $D_n[n+1]$ is to move it to $D_n[1]$ as follows.
For any set $J =\{j_1, j_2, \ldots , j_k\} \subset \{1,2,\ldots , n+1 \}$ and for any integer $i$, let $J+i$ denote the set $\{j_1+i, j_2 +i, \ldots, j_k +i \}$ where the elements are reduced modulo $n+1$ to lie between $1$ and $n+1$ inclusive. Then $\Theta_{n+1}(p_J)=p_{J+1}$, and $(\Theta_{n+1})^i(p_J)= p_{J+i}$.

\section{The link of the vertex $e_{\emptyset}$ and the rotations of the pyramids}

In the $(n+1)$-cube $C_{n+1}$, we have a notion of a {\it link} of any particular vertex.  Here we only  need to examine the concept for the vertex $e_{\emptyset}=(0,0,\ldots, 0)$. The link 
$L(e_{\emptyset})$ consists of those $n$-dimensional cubes that do not contain $e_{\emptyset}$.
These cubes are the intersections of $C_{n+1}$ with the hyperplanes:
 $x_1=1$, $x_2=1$, \ldots $x_n=1$ and $x_{n+1}=1$.  Clearly, $e_{\emptyset}$ is not a point in any of these planes.
 Addresses in each of these planes can be determined by the remaining unspecified vectors.  Thus to find your way
 around $x_1=1$, you  can use the direction vectors, $e_2, e_3,$ through $e_{n+1}$, and in $x_j=1$ you  can use $e_1$,
 through $e_{j-1}$, and $e_{j+1}$ through $e_{n+1}$. Finally in $x_{n+1}=1$ use $e_1$ through $e_n$.

\begin{sloppypar}
The rotation $\Theta_{n+1}$, then rotates these planes among each other. The frame of reference $(e_1, e_2, \ldots, e_n)$ in $x_{n+1}=1$ rotates to  the frame of reference
$(e_2,  e_3, \ldots , e_{n}, e_{n+1})$ in $x_1=1$. Similarly, the frame of  reference $(e_2,e_3, \ldots, e_{n}, e_{n+1})$  for $x_1=1$ transforms to the frame $(e_3, e_4, \ldots, e_n, e_1)$ for $x_2=1$.
\end{sloppypar}

These descriptions of the transformations indicate how the $n$-cubical base of the  pyramid ${\mbox{\rm Py}}_{n+1}$ transforms under the rotations. Since the pyramid is formed by coning this base to $e_{\emptyset}$, this gives
 $(n+1)$-copies of the pyramid. They fill the $(n+1)$-cube as we are coning the link $L(e_{\emptyset})$ 
to $e_{\emptyset}$. Thus ${\mbox{\rm Vol}}({\mbox{\rm Py}}_{n+1})=1/(n+1)$. In even dimensions, the pyramids under go reflections at each stage of the transformation.

This completes the proof that $\int_0^1 x^n \ dx = 1/(n+1)$ for $n>0$. In the next, brief section we discuss how the case $n=0$ fits into  the same framework. In the final section, we indicate how to compute the hyper-volumes of some associated $n+1$-simplices.

\begin{figure}[htb]
\begin{center}
\mbox{
\epsfxsize=4in
\epsfbox{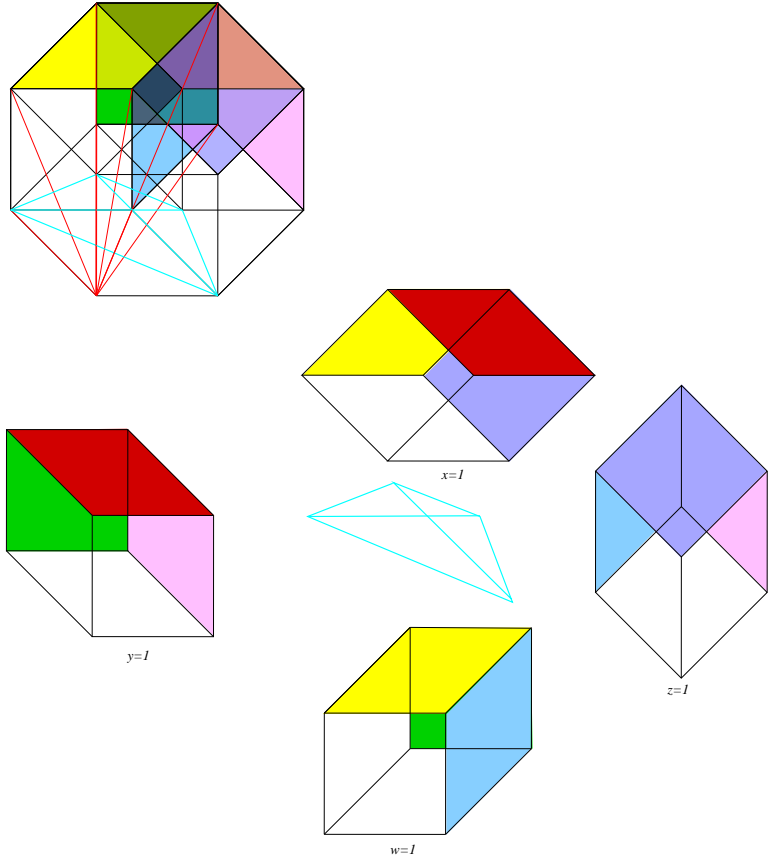}
}
\end{center}
\caption{The link of $(0,0,0,1)$ consists of four cubes of dimension $3$}
\label{linkingcubes}
\end{figure}

\section{Decomposing $C_{n+1}$ into simplices}

Observe that the square $[0,1] \times [0,1]$ is cut into two triangles via the pyramid [sic triangle] ${\mbox{\rm Py}}_{2}$ and its reflection. We can use that decomposition, in the pyramid ${\mbox{\rm Py}}_{3}$ to first decompose the  base into two triangles, and subsequently decompose the pyramid into two tetrahedra. One is the reflection of the other; they are obtained by ``coning off''  the diagonal of the ceiling to the point $e_{\emptyset}=(0,0,0)$.  In this way we can cut the cube into  six congruent pieces, each of  which is a tetrahedron.  The base of the tetrahedron has area $1/2$. the height of each is $1$.  The volume, consequently is $1/6$.

The process continues by induction.  The $(n+1)$-cube can be decomposed into a total  of $(n+1)!$  congruent $(n+1)$-simplices,  each of which has hyper-volume $1/(n+1)!$ The construction follows by decomposing the base cube of the pyramid ${\mbox{\rm Py}}_{n+1}$  into $n!$ simplices, coning within the pyramid, and rotating using $\theta$.

This decomposition is interesting for some historical reasons.  It is impossible to cut a $3$-cube of volume $1$ into finitely many pieces, and reassemble the result into a regular tetrahedron of volume $1$. Such a process is known as a ``scissors congruence.'' The impossibility was demonstrated by Max Dehn who used an invariant defined via the dihedral angles of the cube and tetrahedra.  Scissors congruences give important invariants of $3$-dimensional solids. We remark that AC studies these.

JSC uses decompositions of higher dimensional solids to  identify coincidences among certain cohomology classes. These cohomology theories are used to develop invariants of higher dimensional knots, and the cohomology theories are used to describe and analyze algebraic structures.

Therefore, this paper illustrates the value in ``thinking deeply about  simple things."

\section{Analytic description}
\label{last}

This tesellation of the $(n+1)$-cube has a very pretty analytic description. For each $i=1,\ldots,(n+1)$ let 
$${\rm Py}_{n+1}^i=\{(x_1,\ldots,x_{n+1})\in C_{n+1}|\ x_j \leq x_i\  \forall j=1,\ldots,(n+1), j\neq i \}$$
Tha plane $x_i=x$ intersects ${\rm Py}_{n+1}^i$ in the $n$-cube with the length of each of its edges equal 
to $x$. Hence ${\rm Vol}({\rm Py}_{n+1}^i)=\int_0^1x^ndx$. Also 
$\Theta_{n+1}({\rm Py}_{n+1}^i)={\rm Py}_{n+1}^{i+1}$ for all $i=1,\ldots,n$ and 
$\Theta_{n+1}({\rm Py}_{n+1}^{n+1})={\rm Py}_{n+1}^{1}$ and hence $\{ {\rm Py}_{n+1}^{i}|\ i=1,\ldots,(n+1)\}$ 
tesellate the $(n+1)$-cube showing $\int_0^1 x^n \ dx = 1/(n+1)$ for $n>0$

\section{Mathematica code}

The code below was writing with anearlier version of the paper in mind. In the earlier version the cube of size $x^n$ was located in the plane $x_{n+1}=1-x$ so that the pyramid ${\mbox{\rm Py}_{n+1}}$ had its base on the floor of the $(n+1)$-cube. That turned out to make the rotations of the cube appear to be very complicated. On the other hand, it is not a great idea to modify computer code (even inelegant computer code) that runs. So there are inconguencies between the test and the program. 

\begin{verbatim}





vera[prop_] := If[TrueQ[prop], 1, 0];
val[t_, a_, b_] := vera[a <= t]  vera[t < b] ;
vem[t_, a_, b_] := vera[a <= t]  vera[t <= b] ;
lyle[n_, j_, t_] := Binomial[n, 2] t - j + 1;

booze[n_] := Flatten[Table[{i, j}, {i, 1, n - 1}, {j, i + 1, n}], 1];

wine[j_, n_, t_] := ReplacePart[
      ReplacePart[
        ReplacePart[

          ReplacePart[IdentityMatrix[n],
            Cos[ 2Pi lyle[n, j, t] ], {booze[n][[j]][[1]],
              booze[n][[j]][[1]]}],
          Cos[2Pi lyle[n, j, t]], {booze[n][[j]][[2]],
            booze[n][[j]][[2]]}], -Sin[
            2Pi lyle[n, j, t]], {booze[n][[j]][[2]], booze[n][[j]][[1]]}],
      Sin[ 2Pi lyle[n, j, t]], {booze[n][[j]][[1]], booze[n][[j]][[2]]}];
tire[t_, n_] :=
    Sum[val[t, (j - 1)/Binomial[n, 2], j/Binomial[n, 2]]wine[j, n, t], {j, 1,
          Binomial[n, 2] - 1}] +
      vem[t, (Binomial[n, 2] - 1)/Binomial[n, 2], 1]wine[Binomial[n, 2], n,
          t];


(*This cell establishes a set of  rotations in each of n choose 2 directions,
  one performed after the other.
      The rotations can be appied to any higher dimensional figure *)




ramona[m_]:=Table[1/Sqrt[m],{i,0,m-1}]
(* ramona is used to  give a 3-
    dimensional projection. It may not be the best choice *)


Clear[eh] (*  eh is canadian,
  eh? Now you have to  be careful how you choose to execute the
    next few cells. *)

eh[4]=Transpose[{{0,1},{1,1}/Sqrt[2], {1,0},{1,-1}/Sqrt[2]}] (*
    this value of eh is  good to project the 4-
      cube into the standard octahedral projection.
          It is not much good for any thing else. *)


eh[m_]:={Table[Re[E^(2 Pi j I/m)],{j,0,m-1}],
    Table[Im[E^(2 Pi j I/m)],{j,0,m-1}]}
(* This version of eh gives  a 2dimensional projection to the plane where the \
coordinate axes go to the roots of unity.
      These projections tends to  give the most symmetric pictures. *)


eh[m_]:={Table[Re[E^(2 Pi j I/m)],{j,0,m-1}],
    Table[Im[E^(2 Pi j I/m)],{j,0,m-1}],ramona[m]}
(* This  version of eh projects to 3 dimension,
  it includes some strange distortions. *)



eh[m_]:={IdentityMatrix[m][[1]],IdentityMatrix[m][[2]],
    IdentityMatrix[m][[3]]} (*
    The most obvious choice for eh.
        You can use this to gain more intution since various details are \
squished out.*)

In[8]:=

eh[m_]:={IdentityMatrix[m][[1]],IdentityMatrix[m][[2]],
    IdentityMatrix[m][[3]]+IdentityMatrix[m][[m]]} (*
    Another obvious choice for eh. Try running with this value.
          This is my favorite view for the 4-
      cube. You really can see the cuboid in the simplex here 
even though  other details are squashed out.*)



(* CAUTION. ONLY EXECUTE ONE VALUE OF eh.
      ALSO YOUR CHOICE OF eh DETERMINES WHETHER 
YOU RUN moonbeam or \
moonbeam2d!*)

In[9]:=

topvert[n_]:=Insert[Table[0,{i,1,n-1}],1,n];
missing[j_,m_]:=Join[ Take[Array[x,m],j-1],Take[Array[x,m],j-m]];
nerd[j_,m_,s_]:=Insert[missing[j,m],s,j];
felix[j_,m_,t_]:=
    Line[{Dot[eh[m],tire[t,m]].nerd[j,m,0],
        Dot[eh[m],tire[t,m]].nerd[j,m,1]}];
leonard[j_,m_,s_]:={Hue[.0],
      Line[{Dot[eh[m],tire[t,m]].nerd[j,m,s],
          Dot[eh[m],tire[t,m]].topvert[m]}]};
git[m_]:=Table[IntegerDigits[j,2,m],{j,0,2^m-1}];
shirt[m_]:=Table[IntegerDigits[2j,2,m],{j,0,2^m-1}];
frances[j_,m_,k_,git_]:=
    Table[missing[j,m][[i]]->git[m-1][[k,i]],{i,1,
        Length[missing[j,m]]}];
hyper[j_,m_,t_]:=
    Table[felix[j,m,t]/.frances[j,m,k,git],{k,1,Length[git[m-1]]}];
hypo[j_,m_,t_]:={Hue[j/m^2],
      Table[felix[j,m,t]/.frances[j,m,k,shirt],{k,1,Length[shirt[m-1]]/2}]};
(*This cell provides the main portion of the hyper cube. topvert is e_{n+1}.
        Felix and leonard are lines.
        Most anything with  a feline sounding name is a line: line->
    lion->cat. git stands for digit,
  shirt is the sanitized version of a rhyming word.
        frances is a substitution rule:  rule->
    mule->
      frances the talking mule. It is  misspelled on purpose.
          hypo and hyper are the corresponding cubes.
          The edges of hypo are specially colored. *)



In[19]:=
parrot[j_,k_,m_,
      t_]:={SurfaceColor[
        RGBColor[1-Sqrt[(j^2+k^2)]/(2m)^2,0,1-j/m^2]](*Hue[(j+k)/4m]*),
      Polygon[{Dot[eh[m],tire[t,m]].Array[x,m]/.{x[j]->0,x[k]->0},
          Dot[eh[m],tire[t,m]].Array[x,m]/.{x[j]->1,x[k]->0},
          Dot[eh[m],tire[t,m]].Array[x,m]/.{x[j]->1,x[k]->1},
          Dot[eh[m],tire[t,m]].Array[x,m]/.{x[j]->0,x[k]->1}}]};
twoholes[j_,k_,m_]:=
    Join[Join[ Take[Array[x,m],j-1],Take[Array[x,m],{j+1,k-1}]],
      Take[Array[x,m],k-m]];
fools[j_,k_,m_]:=Table[twoholes[j,k,m][[ell]]->0,{ell,1,m-2}];
rooster[j_,k_,m_,t_]:=parrot[j,k,m,t]/.fools[j,k,m];
hone[j_,m_,t_]:=
  Table[leonard[j,m,0]/.frances[j,m,k,git],{k,1,Length[git[m-1]]}]
(* Things avian are surfaces: It's a bird its a plane. hone is the cone.
          I use the Dr.
          Suess approach to programming to avoid Mathematica thinking that \
uncapitalized words are misspelling of its commands.
          fools are a set of substitution rules.*)

In[24]:=



jag[u_,v_,m_,t_]:={Hue[0.5],
      Line[{Dot[eh[m],tire[t,m]].(IdentityMatrix[m][[u]]+topvert[m]),
          Dot[eh[m],tire[t,m]].(IdentityMatrix[m][[v]]+topvert[m])}]};
uar[u_,m_,t_]:={Hue[0.5],
    Line[{Dot[eh[m],tire[t,m]].(IdentityMatrix[m][[u]]+topvert[m]),
        Dot[eh[m],tire[t,m]].(IdentityMatrix[m][[m]]-topvert[m])}]};
pimplex[m_,t_]:={AbsoluteThickness[3],
    Join[Table[jag[u,v,m,t],{u,1,m},{v,u+1,m-1}],
      Table[uar[u,m,t],{u,1,m-1}]]}
(* jag and uar are sets of lines. pimplex is the  simplex.*)





In[26]:=

powder[m_]:=Table[Reverse[IntegerDigits[k,2,m]] ,{k,2^(m-1)+1,2^m-1}]
(* powder was my approach to  finding the power set.
      This portion of  the program is
    very unelegant. *)



In[27]:=

buns[k_,m_]:=
  Sum[powder[m][[k,ell]],{ell,1,m}] (* buns, sums,
    buns are made of flour a type of powder *)

In[28]:=
mist[m_]:=Do[spot[j,m]={};
      Do[If[Equal[buns[k,m],j],AppendTo[spot[j,m],powder[m][[k]]]],{k,1,
          Length[powder[m]]}],{j,2,m}];
(* mist is a given list.  This is a very weird peice of programming.
      See the next comment *)


In[29]:=
Table[mist[k],{k,3,7}] (*
  OK YOU HAVE TO RUN THIS FOR THE LAST NUMBER AT THE DIMENSION YOU WANT \
TO  GRAPH. THE PROGRAM WORKS IN THEORY FOR ALL N,
  BUT ON MY MACHINES I NEVER GOT PAST 7. I NEED A NEW COMPUTER! IGNORE THE \
NULL OUTPUT.*)

Out[29]=
{Null,Null,Null,Null,Null}

In[30]:=


Clear[freckle]

In[31]:=
freckle[k_,m_,ind_]:= 1/k spot[k,m][[ind]] + (k-2)/k IdentityMatrix[m][[m]]
(*freckle[k,m,ind] is a barycenter of the k-1 dim'l simplex in m-
      space. The index=
    ind is a numbering thereof. 2\[LessEqual]k\[LessEqual]m,
  and ind\[LessEqual]
    Binomial[m-1,k-1].
        The freckle is the spot where the simplex intersects the cone *)


In[32]:=


Clear[victor]

In[33]:=

victor[k_,ell_,m_,ind_,dex_]:=spot[ell,m][[dex]]-spot[k,m][[ind]]
(* What's the vector victor? What's  the clearance, Clarence? ---Airplane *)



In[34]:=

Clear[guess];ClearAll[dorothy]

In[35]:=

guess[m_]:=
  Do[dorothy[k,m]={};
    Do[If [victor[k,k+1,m,ind,dex]\[Equal]Abs[victor[ k,k+1,m,ind,dex]],
        AppendTo[dorothy[k,m],{ind,dex}]],{ind,1,Binomial[m-1,k-1]},{dex,1,
        Binomial[m-1,k]}] ,{k,2,m-1} ]
(* dorothy is supposed to  be close to the cowardly lion.
      dorothy gives the set of lines of intersections between the faces \
of  the cone and the simplex. *)

In[36]:=

Table[guess[m],{m,2,7}] (*
  I guess I  got this part right. It looks  good in low dimensions.
      THE LST NUMBER SHOUlD BE EQUAL TO OR EXCEED THE DIMENSION YOU WANT TO \
GRAPH! BAD PROGRAMMING *)

Out[36]=
{Null,Null,Null,Null,Null,Null}

In[37]:=
Clear[cello]

In[38]:=
cello[k_,ell_,m_,ind_,dex_,t_]:={AbsoluteThickness[2],
      Hue[0.7](*Hue[1-(k+ell)/m^3]*),
      Line[{Dot[eh[m],tire[t,m]].freckle[k,m,ind],
          Dot[eh[m],tire[t,m]].freckle[ell,m,dex]}]};
(* I think I wanted violins here,
  but I  used a cello instead. I  can't remember why.
      You can change the Hue if you like.
      I liked .7 since it was a nice blue.*)

In[39]:=


vegetable[k_,ell_,m_,t_]:=
  Partition[
    Flatten[Table[
        cello[k,ell,m,dorothy[k,m][[j,1]],dorothy[k,m][[j,2]],t],{j,1,
          Length[dorothy[k,m]]}]],3] (*
    vegetables such as carrot sticks. These are the sticks in the simplex*)

In[40]:=
salad[m_,t_]:=
  Table[vegetable[k,k+1,m,t],{k,2,m-1}] (* a bunch of vegetables *)

In[41]:=

polka[k_,m_,ind_,t_]:={PointSize[0.01],Hue[0.99],
      Point[Dot[eh[m],tire[t,m]].freckle[k,m,ind]]};
(* a bunch of dots *)

In[42]:=

dance[m_,t_]:=
  Partition[
    Flatten[ Table[polka[k,m,ind,t],{k,2,m},{ind,1,Length[spot[k,m]]}]],3]
(* such as a polka *)

In[43]:=

ohmy[m_,t_]:= {PointSize[0.02],Point[Dot[eh[m],tire[t,m].Table[0,{j,1,m}]]]}
(* the origin of the cube indicated in a big black dot. *)



(* RUN THIS IF eh  is 2 by m *)
moonbeam2d[m_,dt_,a_,b_]:=
  Table[Show[
      Graphics[{Table[{hyper[j,m,t]},{j,1,m}],{AbsoluteThickness[2],
            Table[{hypo[j,m,t]},{j,1,m-1}]},
          hone[m,m,t],
          salad[m,t],dance[m,t],pimplex[m,t],ohmy[m,t] }],
      PlotLabel->
        StyleForm[IntegerPart[t Binomial[m,2]/dt],"Section",
          FontColor->GrayLevel[.1]],AspectRatio->1,Axes->False,
      ImageSize-> {(40*6),(40*6)}(*,BoxRatios->{1,1,1},
        Boxed->False*)],{t,a,b,dt/Binomial[m,2]}]
(* m is the dimension (n+1)  in the paper.
      dt is an increment of the movies a and b should be between 0  and \
1.  you can control these if you want one rotation *)


moonbeam2d[4,0.0125,0,1]
(* THE LINE THAT GIVES OUTPUT IN CASE eh is 2 by m *)

In[44]:=


(* RUN THIS IF eh  is 3 by m *)
moonbeam[m_,dt_,a_,b_]:=
  Table[Show[
      Graphics3D[{Table[{hyper[j,m,t]},{j,1,m}] ,{AbsoluteThickness[3],
            Table[{hypo[j,m,t]},{j,1,m-1}]},{AbsoluteThickness[2],
            hone[m,m,t]},Table[rooster[j,k,m,t],{j,1,m-1},{k,j+1,m-1}],
          salad[m,t],dance[m,t],pimplex[m,t],ohmy[m,t]}],
      PlotLabel->
        StyleForm[IntegerPart[t Binomial[m,2]/dt],"Section",
          FontColor->GrayLevel[.1]],AspectRatio->1,Axes->False,
      ImageSize-> {(72*6),(72*6)},BoxRatios->{1,1,1},
      Boxed->False],{t,a,b,dt/Binomial[m,2]}]
(* m is the dimension (n+1)  in the paper.
      dt is an increment of the movies a and b should be between 0  and \
1.  you can control these if you want one rotation *)


moonbeam[4,.0125,0,1] (* THE LINE THAT GIVES OUTPUT IN CASE eh is 3 by m *)
\end{verbatim}

\section{A gallery of imagery}

Our inclination in this section is to provide a collection of images that were produced by minor variants on the mathematica code provided. The variations occur in the choices used for the projection matrix {\begin{verbatim} eh. \end{verbatim}} Not every choice that we used are listed in the program code. However, we include some minor comments.

In the first illustration, the figure resembles the drawn image. The reason is that in both cases the octohedral view point was chosen. The rough sketches of the  illustration were done with the computer off, and the Mathematica illustration was  used to check that coordinates were  correctly indicated. The cyan in the drawing was chosen to coincide with the color choice in Mathematica. The next two illustrations are from the same movie.

\begin{center}
\mbox{
\epsfxsize=4in
\epsfbox{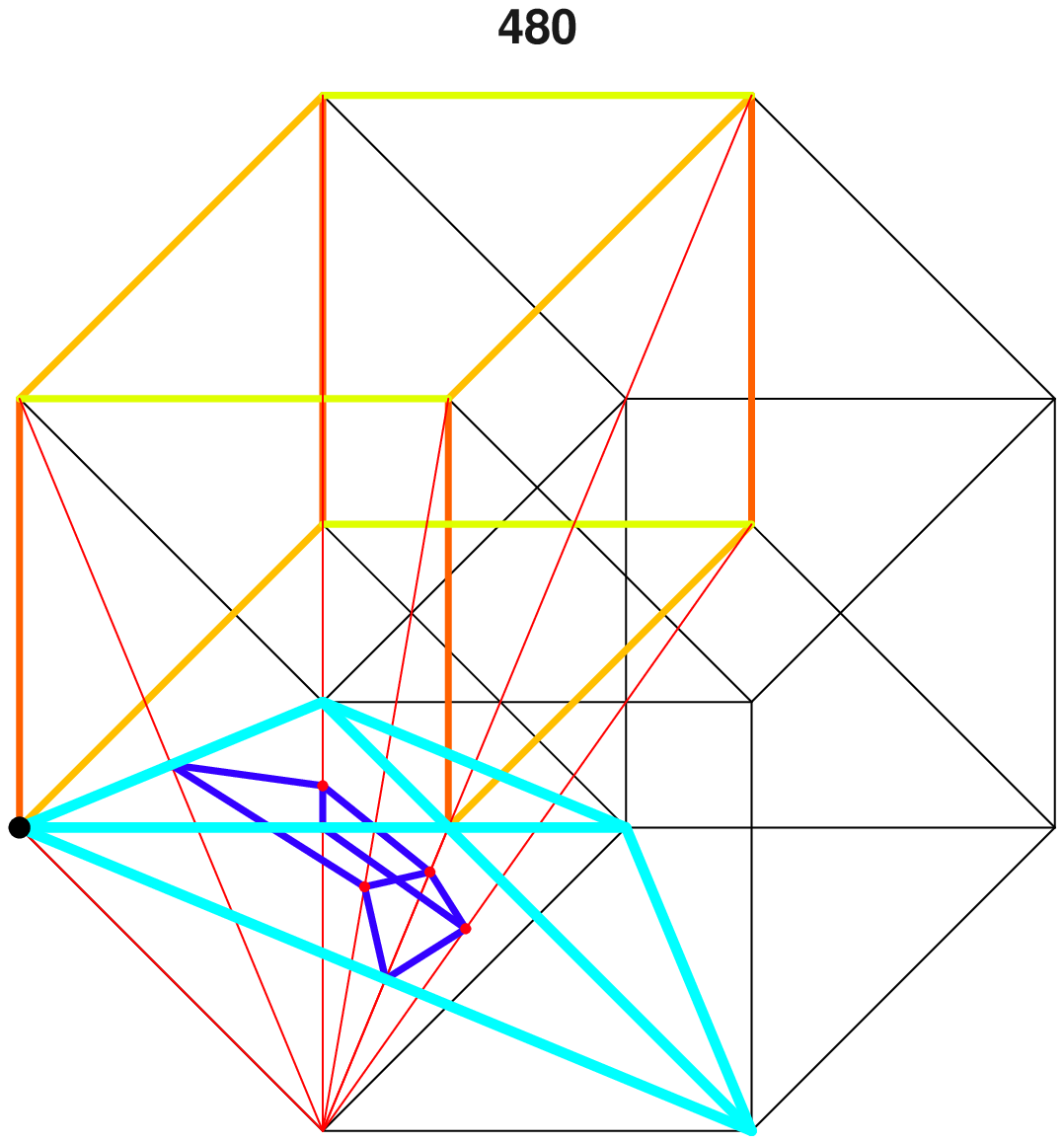}
}
\end{center}

\begin{center}
\mbox{
\epsfxsize=4in
\epsfbox{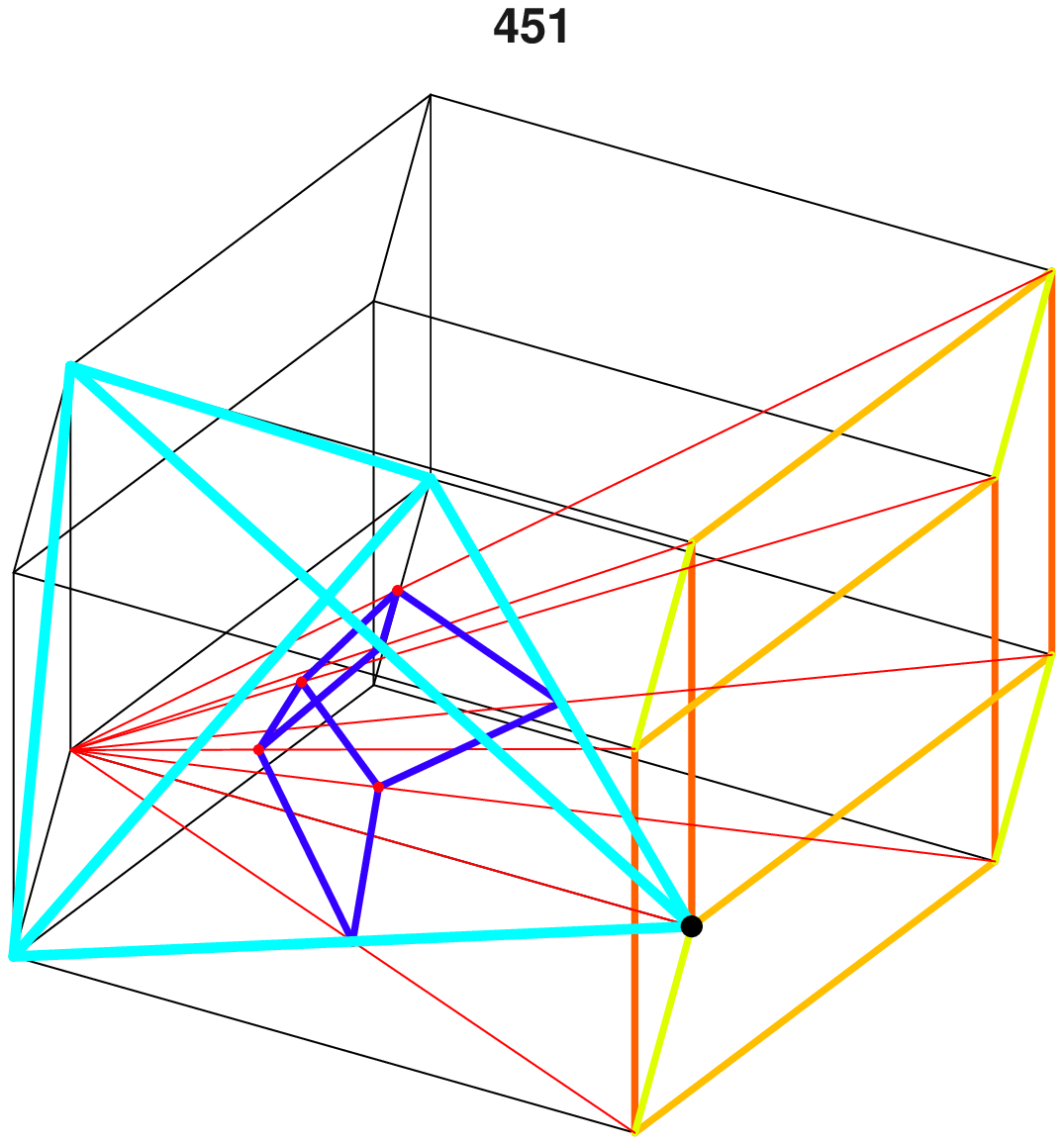}
}
\end{center}

\begin{center}
\mbox{
\epsfxsize=4in
\epsfbox{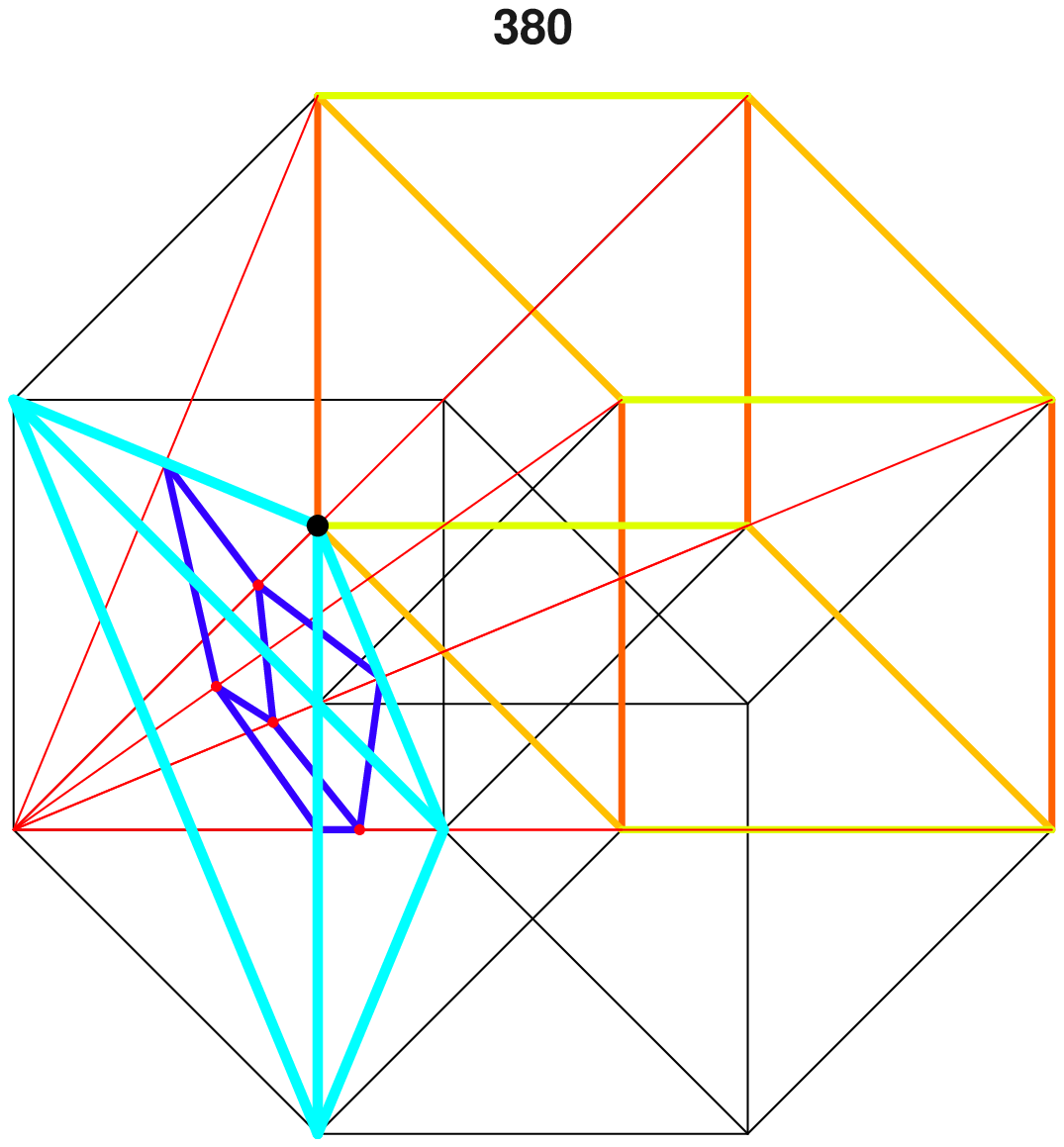}
}
\end{center}

\newpage
The next figure helps the reader see that the intersection of the pyramid with the tetrahdron is a cuboid.

\begin{center}
\mbox{
\epsfxsize=3.5in
\epsfbox{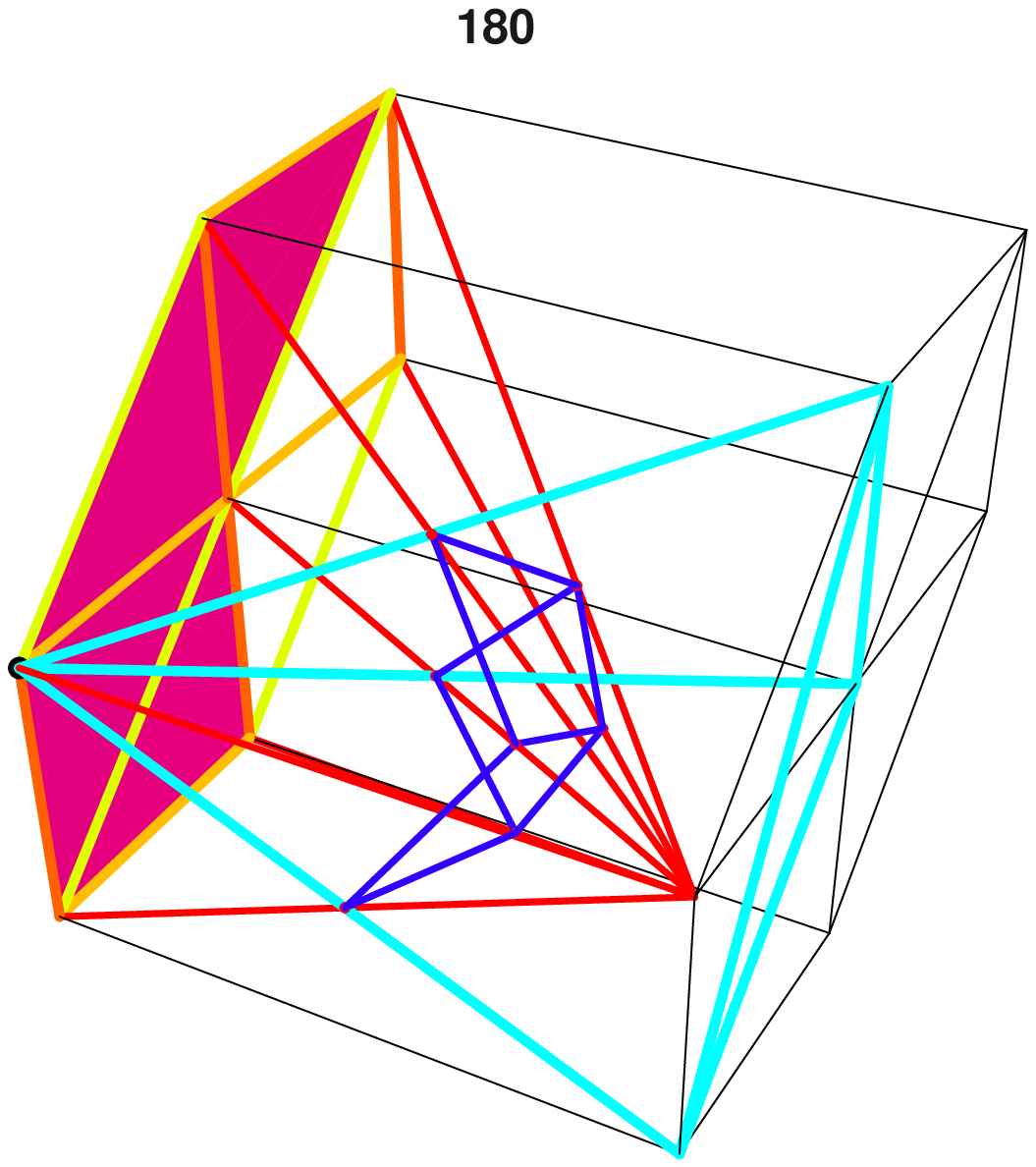}
}
\end{center}

This illustration was chosen using the squished view of the four simplex. The tetrahedron is squashed into a triangle, but the cone to the floor is quite visible.

\begin{center}
\mbox{
\epsfxsize=3.5in
\epsfbox{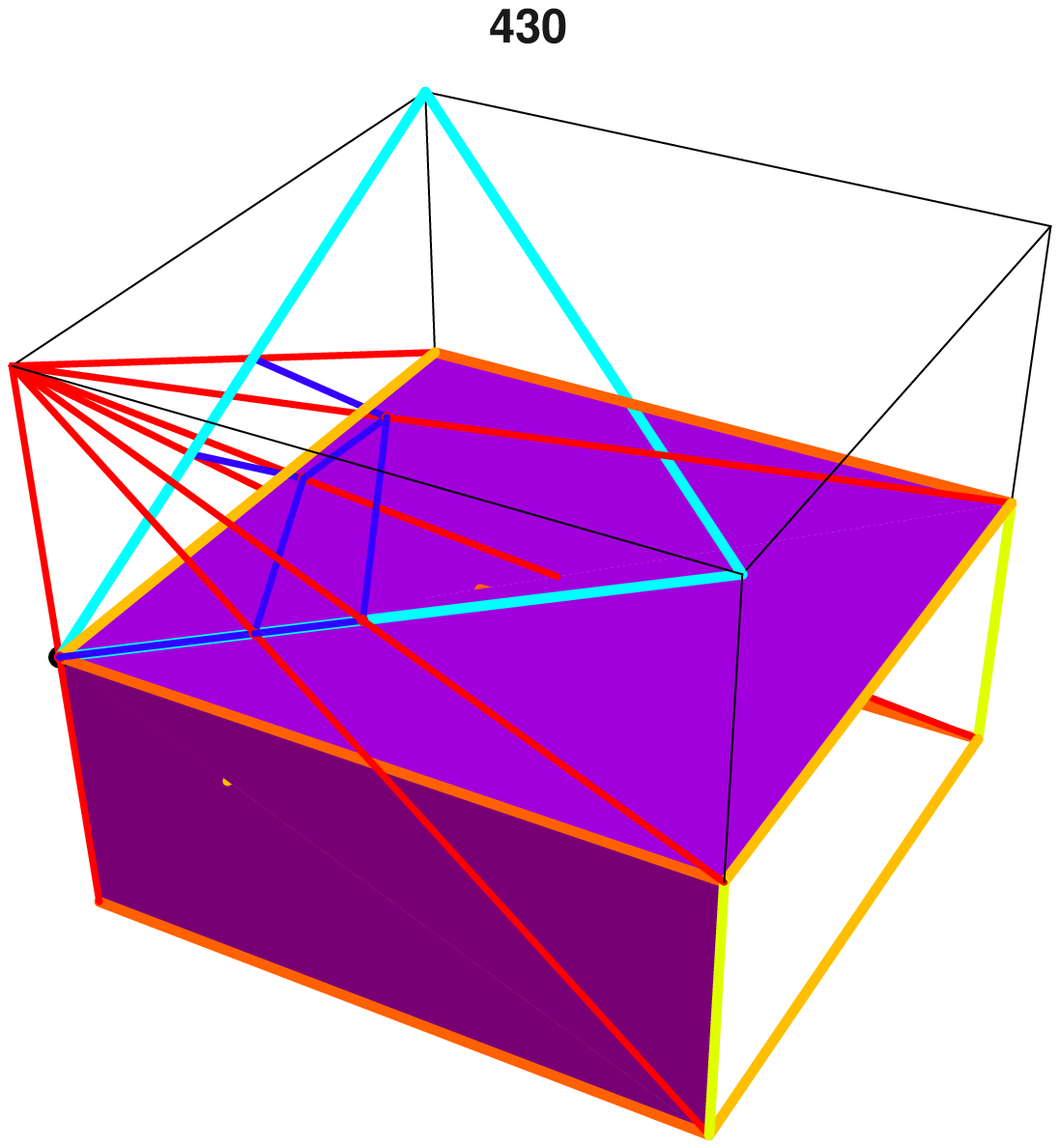}
}
\end{center}

\newpage
Here is one nice view  of the $5$-cube with cone simplex and intersection indicated.

\begin{center}
\mbox{
\epsfxsize=6.5in
\epsfbox{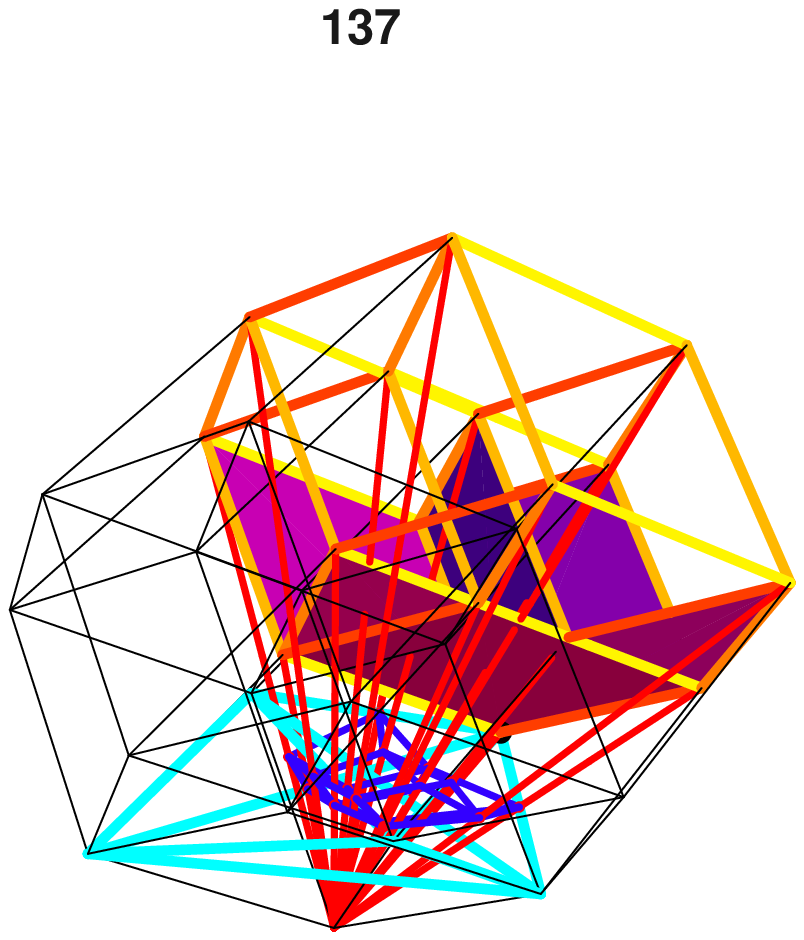}
}
\end{center}

\newpage
This view of the six dimensional cube uses a different value for {\begin{verbatim} eh. \end{verbatim}} On our machines, we can run the corresponding illustration with a large number of stills per movie, but then we cannot animate the result.

\begin{center}
\mbox{
\epsfxsize=6.5in
\epsfbox{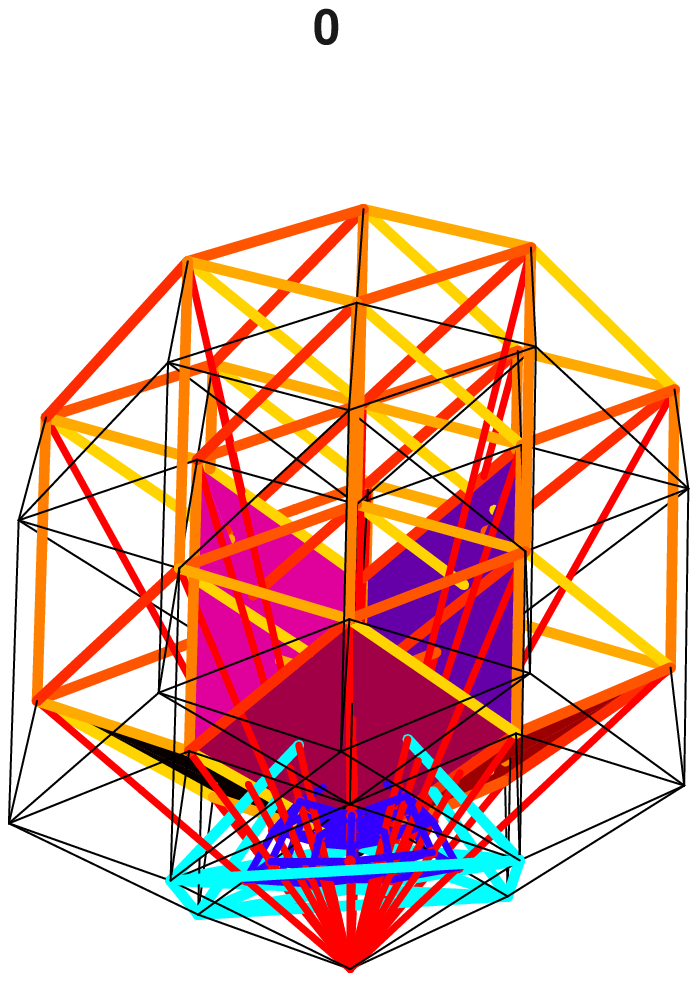}
}
\end{center}

In principle, the  program code listed below runs to illustrate any $n$-cube. In practice it gets bogged down at about $n=6$. The top of the code can be used to illustrate any higher dimensional object. The object is rotated full circle in each of the standard ${n}\choose{2}$ planes in $n$-space. Within the hyper-cubes, only certain $2$-dimensional faces are shown again  ${n}\choose{2}$  of them.

\newpage
Here is a fake stereopsis that was created by taking two successive views of the $4$-cube and setting them side-by-side. Try crossing your eyes to see it in $3D$.

\begin{center}
\mbox{
\epsfxsize=6.5in
\epsfbox{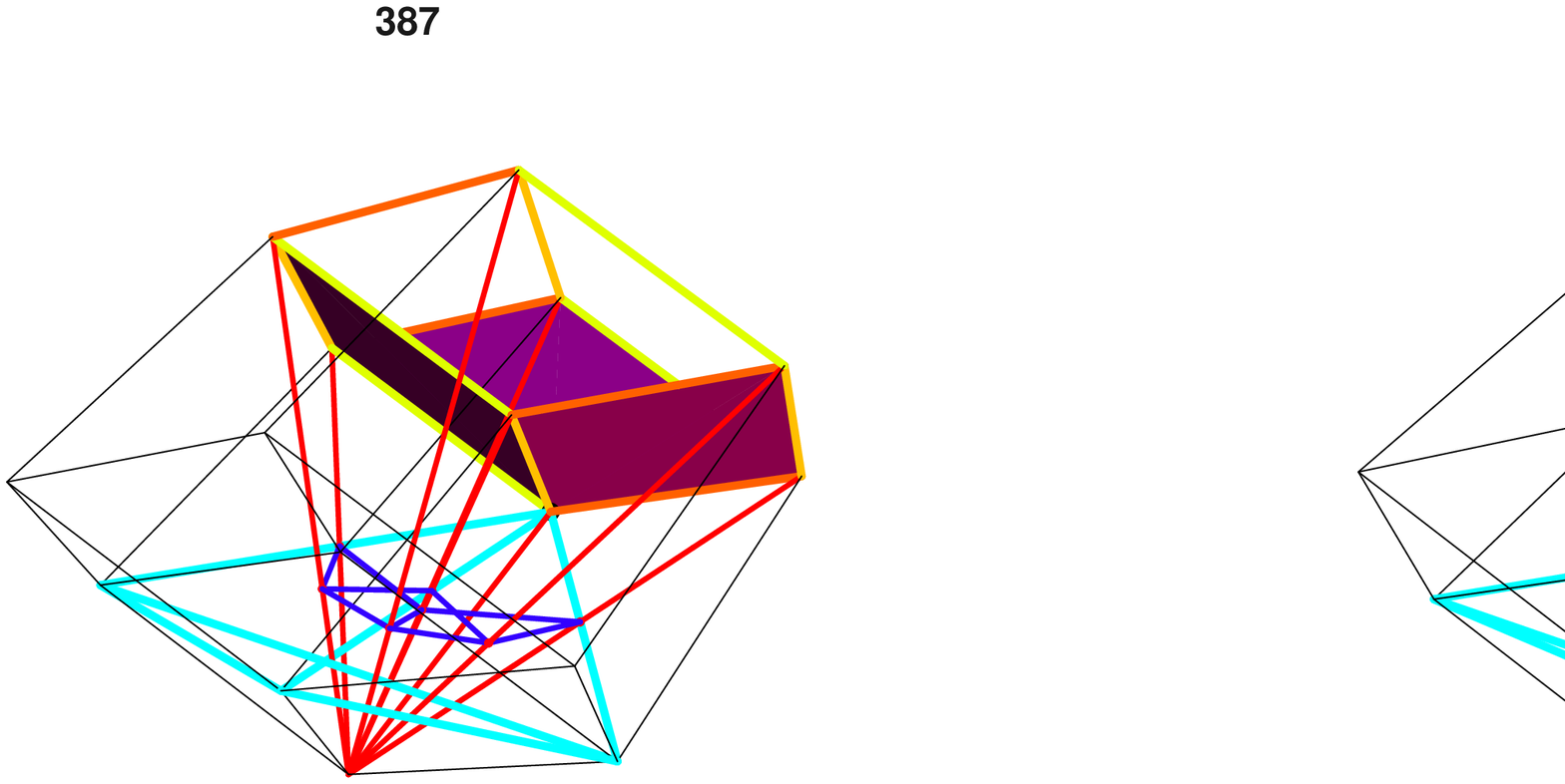}}
\end{center}
 \end{document}